\documentclass{amsart}
\usepackage{amscd}
\newtheorem{theorem}{Theorem}[section]
\newtheorem{lem}[theorem]{Lemma}
\newtheorem{prop}[theorem]{Proposition}

\theoremstyle{definition}

\newtheorem{question}[theorem]{Question}

\newtheorem{claim}[theorem]{Claim}
\newtheorem{conjecture}[theorem]{Conjecture}

\theoremstyle{remark}

\numberwithin{equation}{section}

\newcommand{\Hor}{{\mathcal{H}}}
\newcommand{\V}{{\mathcal{V}}}

\newcommand{\nb}{\nu(\Sigma)}
\newcommand{\nbp}{\nu_p(\Sigma)}

\newcommand{\ra}{\rightarrow}

\newcommand{\Soul}{\Sigma}

\newcommand{\ddt}{\frac{\text{d}}{\text{dt}}}
\newcommand{\ddr}{\frac{\partial}{\partial r}}

\newcommand{\dist}{\text{dist}}

\newcommand{\lb}{\langle}
\newcommand{\rb}{\rangle}

\newcommand{\R}{\mathbb{R}}
\newcommand{\Z}{\mathbb{Z}}

\newcommand{\Ht}{{\hat{\Theta}}}
\newcommand{\RN}{R^{\nabla}}
\begin{document}

\title{Nonnegatively curved metrics on $S^2\times\R^2$}

\author{Detlef Gromoll and Kristopher Tapp}
\address{Department of Mathematics\\ SUNY Stony Brook\\
         Stony Brook, NY 11794-3651}
\email{ktapp@math.sunysb.edu}
\email{detlef@math.sunysb.edu}

\date{\today}
\begin{abstract}
We classify the complete metrics of nonnegative sectional curvature
on $M^2\times\R^2$, where $M^2$ is any compact 2-manifold.
\end{abstract}
\maketitle

\section{Introduction}\label{intro}
According to the Soul Theorem, any open manifold $N^{n+k}$ with nonnegative sectional curvature is diffeomorphic
to the total space of the normal bundle of a compact totally geodesic submanifold $M^n\subset N^{n+k}$~\cite{CG}.
Since this result, some progress has been made on the topological classification problem: for which vector
bundles can the total space admit a complete metric of nonnegative sectional curvature?  However, the metric classification
problem has received little attention beyond the case $n+k\leq 3$, which was completely solved in~\cite{CG}.  We believe that
Perelman's Theorem~\cite{P}, which states that the metric projection $\pi:N\ra M$ is a Riemannian submersion, opens the
door for continued progress on the metric classification.

The simplest remaining case is $n=k=2$, which is most interesting when $M=S^2$.  Every
$\R^2$ bundle over $S^2$ admits nonnegative curvature, and the nontrivial bundles admit very large families of nonnegatively
curved metrics; see~\cite[Section 2]{W} and Section~\ref{rigidity} of this paper.  In contrast, nonnegatively curved metrics
on the trivial bundle $S^2\times\R^2$ are very rigid. Our main theorem is the following classification:

\begin{theorem}\label{main}
The space $M^2\times\R^2$ with an arbitrary complete metric of nonnegative sectional curvature
is isometric to a Riemannian quotient of the form $((M^2,g_0)\times(\R^2,g_F)\times\R)/\R$.
Here $\R$ acts diagonally on the product
by the flow along Killing fields on $(M^2,g_0)$ and $(\R^2,g_F)$ and by translations on $\R$.
\end{theorem}

Our theorem settles the first unknown case of the following open problem:
\begin{question}\label{Qu}
Classify the complete metrics with nonnegative sectional curvature on the trivial vector bundle $S^n\times\R^k$.
\end{question}

Trivial vector bundles represent the simplest case in which one might hope to
completely understand the geometry.  One consequence of Perelman's Theorem is that the boundary of a small metric tube
about a soul is convex, and therefore inherits a metric of nonnegative curvature~\cite{GW}.  Hence, an answer to
question~\ref{Qu} would provide a classification of those nonnegatively curved metrics on $S^n\times S^{k-1}$ which appear
as the boundary of a metric tube about a soul of $S^n\times \R^k$.  We suggest that the following problem is
interesting and more accessible than the full Hopf conjecture:

\begin{conjecture}[Weak Hopf Conjecture] There does not exist a metric of nonnegative curvature on $S^n\times\R^k$
for which the boundary of a small metric tube about the soul has positive curvature in the induced metric.
\end{conjecture}

In Section~\ref{rigidity}, we highlight a precise sense in which any nonnegatively curved metric on $S^2\times\R^2$ is
rigid at the soul, and we ask whether an analogous rigidity must hold for nonnegatively curved metrics
on $S^n\times\R^k$, which is perhaps the correct first step towards resolving the above question and conjecture.

\section{One dimensional Riemannian Foliations}

We begin the proof of Theorem~\ref{main} in this section with some initial lemmas.
A Riemannian submersion with one-dimensional fibers is called homogeneous if the fibers are the integral curves of
a Killing field on the total space.  It was shown in~\cite{GG} that when the total space has constant curvature,
a Riemannian submersion with one-dimensional fibers must be either homogeneous or flat.  We require
the following similar fact, in which $g_0$ denotes an arbitrary metric on $S^2$:

\begin{lem}\label{mainlem}
A Riemannian submersion with one-dimensional fibers and total space $(S^2,g_0)\times\R$ (with the product
metric) must be homogeneous.
\end{lem}

\begin{proof}
Let $\pi:(S^2,g_0)\times\R\ra M^2$ be a Riemannian submersion.
Let $\Hor,\V$ denote the horizontal and vertical spaces of $\pi$.
For $(p,t)\in S^2\times\R$, let $F_{(p,t)}=\pi^{-1}(\pi(p,t))$ denote the $\pi$-fiber through $(p,t)$.

There exists a point $P\in S^2$ such that $\Hor_{(P,0)}=T_{(P,0)}S^2\subset T_{(P,0)} (S^2\times\R)$.
For if there were no such point, then $p\mapsto\Hor_{(p,0)}\cap T_{(p,0)}S^2$ would be a nowhere vanishing
line field on $S^2$, which is not possible.  Notice that $\V_{(P,0)}=\text{span}\{\ddr\}$, where $\ddr$
denotes the unit-length vector field tangent to the $\R$ factor of $S^2\times\R$.  For each geodesic
$\gamma(s)$ though $P$ in $(S^2,g_0)$, the geodesic $(\gamma(s),0)$ in $(S^2,g_0)\times\R$ is horizontal.

Consider the vector field $X$ on $S^2$ defined so that $\V_{(p,0)}=\text{span}\{\ddr+X(p)\}$.
At each point $p\in S^2$ where $X(p)$ is well-defined, it is
tangent to the distance circle about $P$ because $\dist(P,p)=\dist(F_{(P,0)},F_{(p,0)})$.
To establish that $X$ is everywhere well-defined, we must prove that $\V_{(p,0)}$ is never
a subspace of $T_{(p,0)}S^2$.  So let $p\in S^2$.  Let $\gamma(t)$ be a parameterization of the
fiber $F_{(P,0)}$ with $\gamma(0)=(P,0)$ and $\gamma'(0)=\ddr$.  This induces a parameterization,
$\beta(t)$, of the fiber $F_{(p,0)}$ with $\beta(0)=(p,0)$.  More precisely, let $\beta(t)$ be the endpoint
of the lift to $\gamma(t)$ of a fixed minimal geodesic in $M^2$ from $\pi((P,0))$ to $\pi((p,0))$.
If it were the case that $\beta'(0)\in\V_{(p,0)}\subset T_{(p,0)}S^2$, then since $\beta'(0)$ would have to be
tangent to the distance circle in $S^2$ about $P$, we would have $\ddt|_{t=0}\dist(\gamma(t),\beta(t))>0$.
This contradicts the fact that $\dist(\gamma(t),\beta(t))$ is constant, by construction.  Therefore, $X$ is
a well-defined vector field.

For $p\in S^2$, parameterize the fiber through $(p,0)$ as $t\mapsto(f_t(p),t)$,
were $t\mapsto f_t(p)$ is a path in $S^2$ with $f_0(p)=p$.  Notice that $f_t:(S^2,g_0)\ra(S^2,g_0)$
is a smooth family of maps, $f_0=\text{Id}$, and $\frac{\text{D}}{\text{dt}}|_{t=0}f_t(p)=X(p)$.

There must be some point $Q\in S^2$ ($Q\neq P$) where two geodesics through $P$ meet at an angle other than
$180^{\circ}$, so that $\V_{(Q,0)}=\text{span}\{\ddr\}$ as well.  Further, we can assume that
there is a unique such point.  For, suppose that $Q_1$ and $Q_2$ were two such points.  Then at any point $y\in S^2$
for which the three vectors in $T_y S^2$ tangent to the minimal paths to $P,Q_1$ and $Q_2$ span $T_yS^2$,
we have that $\V_{(y,0)}=\text{span}\{\ddr\}$.  But this property is true for all $y\in S^2$ outside of
a set of measure zero, so $\V_{(y,0)}=\text{span}\{\ddr\}$ for all $y\in S^2$.
Therefore, for any pair $x,y\in S^2$, $\dist(F_{(x,0)},F_{(y,0)})=\dist(x,y)$.  This implies that for all $t\in\R$,
$\dist(f_t(x),f_t(y))\geq\dist(x,y)$.  Since $S^2$ is compact, $f_t$ must be an isometry for each $t$.
So, for any $x,y\in S^2$ and any $t\in\R$, the minimal path in $(S^2,g_0)$ from $f_t(x)$ to $f_t(y)$
realizes the distance between the fibers $F_{(f_t(x),t)}=F_{(x,0)}$ and $F_{(f_t(y),t)}=F_{(y,0)}$, and is
therefore horizontal.  This means that $\V_{(y,t)}=\text{span}\{\ddr\}$ for all $y\in S^2$ and all $t\in\R$,
so the submersion is trivial.

So we can assume that there is a single point $Q$ such that all geodesics from $P$ meet at $Q$.
By a standard variation argument, any two geodesics from $P$ to $Q$ share a common length.
We argue next that $f_{t_0}(P)=P$ and $f_{t_0}(Q)=Q$ for all $t_0\in\R$.
Since $\dist(P,Q)=\dist(F_{(P,0)},F_{(Q,0)})$, we know that $\dist(f_{t_0}(P),f_{t_0}(Q))\geq\dist(P,Q)$.
Assume that $\dist(P,f_{t_0}(P))\leq\dist(Q,f_{t_0}(Q))$ (the other case is similar).  The minimal path from
$f_{t_0}(P)$ to $P$ followed by the minimal path from $P$ to $f_{t_0}(Q)$ has distance $<\dist(P,Q)$ unless
$\dist(f_{t_0}(P),P)=\dist(f_{t_0}(Q),Q)$ and $\{P,f_{t_0}(Q),Q,f_{t_0}(P)\}$ all lie on a single geodesic loop.  In this
case, $\dist(f_{t_0}(P),f_{t_0}(Q))=\dist(P,Q)$, and $\frac{\text{D}}{\text{dt}}|_{t=t_0} f_t(P)$ must be tangent to a
distance circle about $P$; for otherwise we would have that $\ddt|_{t=t_0}\dist(f_t(P),f_{t_0}(Q))<0$, and so
$\ddt|_{t=t_0}\dist((f_t(P),t),(f_{t_0}(Q),t_0))<0$, which would put the fibers $F_{(P,0)}$ and $F_{(Q,0)}$ too close
together.  So $\frac{\text{D}}{\text{dt}}f_t(P)$ is everywhere tangent to the distance circles about $P$, which implies that
$f_t(P)=P$ (and similarly $f_t(Q)=Q$) for all $t\in\R$. In particular, this implies that each path $f_t(p)$ is contained in a
single distance circle about $P$.

We prove next that $X$ is a Jacobi field along any geodesic $\alpha(s)$ in $S^2$
through $P$.  Notice that $X(\alpha(s))$ is the variational vector field of the family
$(t,s)\mapsto f_t(\alpha(s))$, which we wish to show is a family of geodesics.  For this,
it will suffice to prove that $\dist(f_t(\alpha(s_1)),f_t(\alpha(s_2)))=\dist(\alpha(s_1),\alpha(s_2))$.
But if this were not the case, then we would have:
$$\dist((f_t(\alpha(s_1)),t),F_{(\alpha(s_2),0)}) > \dist(\alpha(s_1),\alpha(s_2))
  =\dist(F_{(\alpha(s_1),0)},F_{(\alpha(s_2),0)}),$$
which would provide a contradiction.  Therefore $X$ is a Jacobi field along any geodesic through $P$,
which implies that $X$ is the exponential image of a rotational vector field on $T_P S^2$.

For fixed $t\in\R$, define the vector field $X_t$ on $S^2$ so that $\V_{(p,t)}=\text{span}\{\ddr+X_t(p)\}$.
By the above arguments, each $X_t$ is the exponential image of a rotational vector field on $T_P S^2$
(also on $T_Q S^2$); hence, each $X_t$ is a multiple of $X=X_0$.  Let $A$ denote a single distance circle about
$P$, and let $C=\{(p,t)\in S^2\times\R\mid p\in A\}=A\times\R$.  Notice that $\pi|_C$ is a Riemannian submersion
from $C$ onto some circle in the base space $M^2$.  Since $C$ is flat, the fibers of this submersion must be tangent
to a Killing field on $C$; see~\cite{GG} in which one dimensional Riemannian foliations of constant curvature spaces
are classified.  This proves that (1) $X$ has constant norm along $C$, and (2) $X_t=X$ for all $t\in\R$.
Since $X$ has constant norm along each circle $A$, it is a Killing Field.  This completes the proof.
\end{proof}

Theorem~\ref{main} concludes that $M^2\times\R^2$ with an arbitrary metric of nonnegative curvature
is isometric to a Riemannian quotient of the form $((M^2,g_0)\times(\R^2,g_F)\times\R)/\R$.  So, the distance spheres
about the soul must themselves be quotient metrics of the form
$((M^2,g_0)\times S^1(r)\times\R)/\R$, where $S^1(r)$ denotes the circle with circumference
$2\pi r$.  In the remainder of this section, we do some calculations to see metrically what
these distance spheres and their projections onto the soul look like.  These calculations
will be used in the proof of Theorem~\ref{main}

More precisely, consider the following $\R$ action on $(M,g_0)\times S^1(r)\times\R$:
$$(p,\theta,t)\stackrel{s}{\mapsto}(\rho_sp,\theta+s,t+s),$$
where $\rho_s$ denotes the flow along a Killing field $X$ on $(M,g_0)$ for time $s$.
This gives the following diagram:
$$(M\times S^1,h)\stackrel{d}{\leftarrow}((M,g_0)\times S^1(r)\times\R)/\R
  \stackrel{\Pi}{\ra} ((M,g_0)\times\R)/\R
  \stackrel{\delta}{\ra} (M,g_{\Soul}).$$
Here the diffeomorphism $d$ is defined as $[p,\theta,t]\stackrel{d}{\mapsto}(\rho_{-t}p,\theta-t)$,
and $h$ is defined as the metric on $M\times S^1$ which makes $d$ and isometry.
The diffeomorphism $\delta$ is defined as $[p,t]\stackrel{\delta}{\mapsto} \rho_{-t}p$,
and $g_{\Soul}$ is defined as the metric on $M$ which makes $\delta$ an isometry.  Finally, $\Pi$
is defined as $[p,\theta,t]\stackrel{\Pi}{\ra}[p,t]$, and is a Riemannian submersion.

The circle bundle map $\phi=\delta\circ\Pi\circ d^{-1}:(M\times S^1,h)\ra (M,g_{\Soul})$
is a Riemannian submersion which has the simple description $(p,\theta)\stackrel{\phi}{\ra}p$.
We now study the metric properties of this circle bundle.

\begin{claim}\label{len2}
$\text{length}(\phi^{-1}(p))=\frac{2\pi r}{\sqrt{1+r^2(1-|X(p)|^2_{g_{\Soul}})}}$
\end{claim}

\begin{proof}[proof of Claim~\ref{len2}]
By~\cite{C}, the metric $h$ is obtained from the product metric $h_0$
on $(M,g_0)\times S^1(r)$ by rescaling along the Killing field $X+\Ht$,
where $\Ht$ denotes the vector field tangent to $S^1$ corresponding to a rotation
at one radian per unit speed.  So $\Ht$ can be decomposed into vectors parallel and
perpendicular to this Killing field as follows:
$$\Ht=a(X+\Ht) + (-aX + (1-a)\Ht),$$
Where $a$ is chosen so that the two vectors, $a(X+\Ht)$ and $-aX+(1-a)\Ht$, are
$h_0$-perpendicular, namely, $a=\frac{r^2}{r^2+|X|_{g_0}^2}$.
Thus,
\begin{eqnarray*}
|\Ht|_{h}^2  & = & \frac{|aX+a\Ht|^2_{h_0}}{1+|X+\Ht|^2_{h_0}}
           +|-aX+(1-a)\Ht|^2_{h_0}\\
         & = & \frac{a^2|X|^2_{g_0} + (ar)^2}{1+|X|^2_{g_0}+r^2}
            + a^2|X|^2_{g_0} + (1-a)^2r^2\\
         & = & \frac{r^2+r^2|X|^2_{g_0}}{1+|X|^2_{g_0}+r^2}\\
         & = & \frac{r^2}{1+r^2(1-|X|^2_{g_{\Soul}})}.
\end{eqnarray*}
The fourth equality follows from
the relationship $|X|^2_{g_{\Soul}} = \frac{|X|_{g_0}^2}{1+|X|_{g_0}^2}$, and the claim
follows because $\text{length}(\phi^{-1}(p))=2\pi|\Ht|_h$.
\end{proof}

\begin{claim}\label{hor2}
The horizontal space of $\phi$ is $\text{span}\{Y,X+|X|^2_{g_{\Soul}}\Ht\}$,
where $Y$ is a (local) vector field on $M$ orthogonal to $X$.
\end{claim}

\begin{proof} [proof of Claim~\ref{hor2}]
It is easy to see that $Y$ is $h$-orthogonal to $\Ht$, and hence orthogonal to the fibers.
A second vector orthogonal to $\Ht$ will have the form $X+c\Ht$, and we show next
that $c=|X|^2_{g_{\Soul}}$ is the correct choice for $c$.  To see this, it is necessary
to decompose both $\Ht$ and $X+c\Ht$ into components parallel and perpendicular
to the Killing field $X+\Ht$:

\begin{eqnarray*}
\Ht & = & [a(X+\Ht)] + [-aX+(1-a)\Ht]\\
(X+c\Ht) & = & [b(X+\Ht)] + [(1-b)X + (c-b)\Ht]
\end{eqnarray*}

As in the proof of the previous lemma, $a=\frac{r^2}{r^2+|X|^2_{g_0}}$.  Similarly, solving for $b$ so that
the two vectors $b(X+\Ht)$ and $(1-b)X + (c-b)\Ht$ are $h_0$-perpendicular gives
$b=\frac{|X|^2_{g_0}+cr^2}{|X|_{g_0}^2+r^2}$.

This allows us to compute:
\begin{eqnarray*}
\lefteqn{\lb \Ht,X+c\Ht\rb_h}\\
 & = & \lb [a(X+\Ht)] + [-aX+(1-a)\Ht],[b(X+\Ht)] + [(1-b)X + (c-b)\Ht]\rb_h\\
 & = & \lb [a(X+\Ht)],[b(X+\Ht)]\rb_h + \lb [-aX+(1-a)\Ht],[(1-b)X + (c-b)\Ht]\rb_{h_0}\\
 & = & \frac{ab(|X|_{g_0}^2 + r^2)}{1+|X|_{g_0}^2+r^2} -a(1-b)|X|^2_{g_0} +(c-b)(1-a)r^2=0
\end{eqnarray*}
Solving for $c$ and simplifying gives:
$$c=\frac{|X|^2_{g_0}}{|X|^2_{g_0}+1} = |X|^2_{g_{\Soul}}.$$
Notice in particular that $c$ does not depend on $r$.
\end{proof}


\section{Proof of Main Theorem}
\begin{proof}[proof of Theorem~\ref{main}]
Let $E$ denote $M^2\times\R^2$ together with some metric of nonnegative sectional curvature.
Let $\Soul\subset E$ denote a soul of $E$, and let $\pi:E\ra\Soul$ denote the metric
projection.

If $M^2$ is a torus, then the soul is flat, so the $A$-tensor of $\pi$ vanishes, which implies
that the metric is locally (and hence also globally) a Riemannian product~\cite[Theorem 1.3]{W2}.  On the other hand,
the theorem when $M^2$ is the $2$-sphere implies the theorem when $M^2$ is $\mathbb{RP}^2$.
We can therefore assume that $M^2$ is the $2$-sphere, so $\Sigma = (S^2,g_{\Sigma})$, where
$\Sigma$ denotes the metric on the soul.

If the normal bundle of $\Soul$ in $E$ is flat, then $E$ is a Riemannian product, so
we can assume that it is not flat.  In this case, the holonomy group is transitive
and the normal exponential map is a diffeomorphism.

Let $S$ denote the boundary of a ball of some radius about $\Soul$ in $E$.  The intrinsic
metric, $h$, on $S$ is smooth and has nonnegative curvature by the Gauss equation.  The
universal cover $\tilde{S}$ of $S$ is diffeomorphic to $S^2\times\R$.
By the splitting theorem, $\tilde{S}$ with its pull-back metric is isometric to
the Riemannian product $(S^2,g_0)\times\R$ for some metric $g_0$ of nonnegative curvature
on $S^2$.  So we have:
$$(S^2,g_0)\times\R\stackrel{f}{\ra}(S^2\times S^1,h)\stackrel{\pi}{\ra}(S^2,g_{\Soul}).$$
Lemma~\ref{mainlem} applied to the
composition $\pi\circ f$ implies that there is a Killing field $X$ on $(S^2,g_0)$, and
that $(S^2,g_{\Soul})$ is just the Riemannian quotient $((S^2,g_0)\times\R)/\R$.
The $\R$ action on $(S^2,g_0)\times\R$ can be described as
$$(p,t)\stackrel{s}{\ra}(\rho_sp,t+as),$$
where $\rho_s$ denotes the flow for time $s$ along $X$, and $a>0$ is some constant.

On the other hand, $f$ is a covering map, so there is an isometric $\Z$-action on
$(S^2,g_0)\times\R$ so that $(S^2\times S^1,h)=((S^2,g_0)\times\R)/\Z$.  The $\Z$ action comes from some
imbedding $\Z\subset\R$, and we can assume without loss of generality that it is the natural embedding, so
for $n\in\Z$,
$$(p,t)\stackrel{n}{\ra}(\rho_np,t+an).$$
A diffeomorphism $d:(S^2\times\R)/\Z\ra S^2\times S^1$ is defined as follows:
$$[p,t]\stackrel{d}{\ra}(\rho_{-t/a}p,[t/a]).$$
Here $[t/a]\in S^1=\R/\Z$.  It is straightforward to show that $d$ is well-defined,
injective and surjective, and that $d^{-1}(p,[t])=[\rho_{t}p,at]$.  We can consider
$h$ to be the metric on $S^2\times S^1$ which makes
$d:((S^2,g_0)\times\R)/\Z\ra(S^2\times S^1,h)$ an isometry.
Additionally, the diffeomorphism $\delta:(S^2\times\R)/\R\ra S^2$ defined so that
$[p,t]\stackrel{\delta}{\ra}\rho_{-t/a}p$ provides an isometry
$\delta:((S^2,g_0)\times\R)/\R)\ra(S^2,g_{\Soul})$.  The story so far is summarized by
the following diagram:

$$
\begin{CD}
  @. (S^2,g_0)\times\R @. @.\\
  @. @VVjV @. @.\\
  (S^2\times S^1,h) @<d<< ((S^2,g_0)\times\R)/\Z @>\Pi>> ((S^2,g_0)\times\R)/\R @>{\delta}>> (S^2,g_{\Soul})
\end{CD}
$$

Notice that the circle bundle map
$\phi=\delta\circ\Pi\circ d^{-1}:(S^2\times S^1,h)\ra (S^2,g_{\Soul})$
has the simple description $\phi((p,[t]))=p$.  Next we study the metric properties
of this circle bundle.

\begin{claim} \label{len1}
  $\text{length}(\phi^{-1}(p))= (\frac{1}{a^2}-\frac{1}{a^4}|X|_{g_{\Soul}}^2)^{-1/2}$.
\end{claim}

\begin{proof}[proof of Claim\ref{len1}]
Notice that the fiber is parameterized as
$\gamma(t)=(p,[t])$, $t\in[0,1]$.  One $(d\circ j)$-lift of $\gamma$ to
$(S^2,g_0)\times\R$ looks like $\tilde{\gamma}(t)=(\rho_tp,at)$, so
$|\tilde{\gamma}'(t)|=\sqrt{a^2+|X(p)|_{g_0}^2}$.  But the metric
$g_{\Soul}$ on $S^2$ is obtained from the metric $g_0$ by rescalling along the
Killing field $X$.  More precisely,
$|X(p)|^2_{g_{\Soul}} = \frac{a^2|X(p)|_{g_0}^2}{a^2+|X(p)|_{g_0}^2}$.
This substitution proves the claim.
\end{proof}

\begin{claim}\label{hor1}
The horizontal space of $\phi$ is:
$$\Hor=\text{span}\{Y,-X+\frac{2\pi}{a^2}|X|^2_{g_{\Soul}}\hat{\Theta}\},$$
where $Y$ is a (local) vector field on $S^2$ perpendicular to $X$, and $\hat{\Theta}$ is the vector
field on $S^1$ corresponding to rotation at 1 radian per unit time.
\end{claim}

\begin{proof}[proof of claim~\ref{hor1}]
Consider the curve $\tilde{\gamma}(t)=(\rho_tp,at)$ on $(S^2,g_0)\times\R$
as in the previous claim.  $\tilde{\gamma}'(0)=X+a\ddr$, where $\ddr$ denotes the unit-length
vector field tangent to the second factor of $(S^2,g_0)\times\R$.  So
$(\tilde{\gamma}'(0))^{\perp} = \text{span}\{Y, -aX + |X|^2_{g_0}\ddr\}$.
For the first vector, $Y$, it's clear that $(d\circ j)_*(Y)=Y$.  For the second vector,
consider the path $\tilde{\alpha}(t)=(\rho_{-at}p,|X|^2_{g_0}t)$, which is defined so that
$\tilde{\alpha}'(0)=-aX+|X|^2_{g_0}\ddr$.
Let $\alpha=d \circ j\circ\tilde{\alpha}$ be the projection of this curve.  But
$$\alpha(t)=(\rho_{\sigma} p,\,(t/a)[|X|^2_{g_0}]),\text{ where } \sigma=-(t/a)(|X|^2_{g_0}+a^2)$$
so $\alpha'(0)=-\frac{|X|^2_{g_0}+a^2}{a}X + 2\pi\frac{|X|^2_{g_0}}{a}\hat{\Theta}$,
which is parallel to $-X+2\pi(\frac{|X|_{g_0}^2}{a^2+|X|_{g_0}^2})\hat{\Theta}$.
Re-writing this vector in terms of $g_{\Soul}$ completes the claim.
\end{proof}

Notice that the metric $h$ on $S^2\times S^1$ is completely determined by
the data $\{g_{\Soul},X,a\}$.  More precisely, $h$ is determined by requiring
that the projection $\phi:(S^2\times S^1,h)\ra(S^2,g_{\Soul})$
onto the first factor is a Riemannian submersion with horizontal space
as in Claim~\ref{hor1} and vertical space $\Ht$ uniformly rescaled so that the fiber lengths
are as in Claim~\ref{len1}.

By Claims~\ref{len2},\,\ref{hor2},\,\ref{len1}, and \ref{hor1}, this metric $h$
on $S^2\times S^1$ can be re-described as a
Riemannian quotient of the form $((S^2,g_1)\times S^1(r)\times\R)/\R$,
where $\R$ acts on $(S^2,g_1)$ by flow along a Killing vector field $\hat{X}$ and
on $S^1(r)$ by rotation.
For this quotient metric, $\Hor=\text{span}\{Y,\hat{X}+|\hat{X}|^2_{g_{\Soul}}\Ht\}$ and
$\text{fiberlength}=(\frac{1+r^2}{4\pi^2 r^2}-\frac{|\hat{X}|^2_{g_{\Soul}}}{4\pi^2})^{-1/2}$,
so it agrees with the metric $h$ if we choose $\hat{X}=-\frac{2\pi}{a^2}X$, and choose $r$ such that
$\frac{1}{a^2}=\frac{1+r^2}{4\pi^2 r^2}$.

So each distance sphere $S$ in $E$ about $\Soul$ can be described as a Riemannian quotient
of the form $((S^2,g_1)\times S^1(r)\times\R)/\R$.  Further, two distance spheres (at different
distances from $\Soul$) in $E$, when re-described as submersion metrics, have the same
soul metric $g_{\Soul}$ and also the same horizontal spaces (because the distribution
determined by a connection is scale invariant).
This implies that there is a single metric $g_1$ and a single Killing field
$\hat{X}$ on $(S^2,g_1)$ such that all distance spheres can be described as
$((S^2,g_1)\times S^1(r)\times\R)/\R$ for varying $r$.  It follows that the whole manifold
$E$ can be described as $((S^2,g_1)\times(\R^2,g_F)\times\R)/\R$, where $g_F$ is a
rotationally-symmetric metric whose distance spheres have lengths prescribed by the
above requirements.
\end{proof}

\section{Rigidity at the soul}\label{rigidity}
In this section we exhibit a sense in which any complete metric of nonnegative curvature on
$S^2\times\R^2$ is rigid at its soul.

Let $M^4$  be an open manifold of nonnegative curvature with soul $\Sigma^2$.
Let $k:\Sigma\ra\R$ denote the sectional curvature function of the soul.  The curvature $\RN$ of the connection $\nabla$ in
normal bundle $\nb$ of $\Soul$ is determined by a function $f:\Soul\ra\R$ (so that for oriented orthonormal bases
$\{X,Y\}$ of $T_p\Soul$ and $\{W,V\}$ of $\nbp$, we have $f(p)=\lb\RN(X,Y)W,V\rb$). The curvature of 2-planes normal to the soul
are determined by a function $g:\Soul\ra\R$ (so that $g(p)=\lb R(W,V)V,W\rb$).  Notice that
$\int_{S^2} f \text{dvol} = 0$ because the bundle is trivial.  It follows from Theorem A
of~\cite{T} that the following inequality is satisfied for all $p\in\Soul$ and all unit-length vectors $X\in T_p\Soul$:
\begin{equation}\label{tapp}
(Xf)^2 \leq (f^2 + \frac{2}{3}\text{hess}_g(X,X))\cdot k.
\end{equation}
This inequality comes from the fact that any mixed
2-plane $\sigma$ at the soul is flat, so the hessian of the sectional curvature function on the Grassmanian of all 2-planes
of $M$ must be nonnegative definite at $\sigma$.

For the natural metric of nonnegative curvature on a nontrivial $\R^2$ bundle over $S^2$,
inequality~\ref{tapp} is strict.  The strictness of the inequality
means two things.  First, the metric is non-rigid.  For example, the metric on the soul and the connection in the normal
bundle of the soul can both be arbitrarily perturbed without losing nonnegative curvature; see~\cite[Theorem B]{T}
or~\cite[Proposition 3.3]{SW}.  Second, the boundary of a small distance sphere
about the soul (which is diffeomorphic to a lens space) has strictly positive curvature in the
induced metric~\cite[Theorem C]{T}.

If, on the other hand, $M^4=S^2\times\R^2$, then the metric is much more rigid.  Inequality~\ref{tapp} could not possibly be strict
for all $X$ because the sphere bundle $S^2\times S^1$ does not admit positive curvature. In fact, it is not strict for any $X$:
\begin{prop}\label{equality}
Any metric of nonnegative curvature on $S^2\times\R^2$ is rigid in the sense that
$(Xf)^2 = (f^2 + \frac{2}{3}\text{hess}_g(X,X))\cdot k$ for any unit-vector $X$ tangent to any point $p$ of the soul.
\end{prop}
\begin{proof}[Idea of proof]
By Theorem~\ref{main}, an arbitrary metric $h$ of nonnegative curvature on $S^2\times\R^2$ can be written as
$$(S^2\times\R^2,h)=((S^2,g_0)\times(\R^2,g_F)\times\R)/\R.$$
The metric $g_0$ is rotationally-symmetric, so it can be described as
a surface of revolution in $\R^3$ parameterized as
$(s,t)\mapsto(f(t)\cdot\cos(s),f(t)\cdot\sin(s),\phi(t))$,
where $\phi'(t)^2 + f'(t)^2 = 1$. The $\R$-action on $(S^2,g_0)$ is determined by some
multiple, $C_1$, of the Killing field which rotates at 1 radian per unit time.
Similarly, the $\R$-action on $(\R^2,g_F)$ is determined by some multiple, $C_2$, of the Killing
field which rotates at 1 radian per unit time.  Let $\Lambda$ denote the sectional curvature
of $(\R,g_F)$ at the vertex.

It is then straightforward to compute each term of inequality~\ref{tapp}
in terms only of $\{f,\Lambda,C_1,C_2\}$.  The proposition then follows by direct calculation.
\end{proof}

It is somewhat unsatisfying that we used Theorem~\ref{main} to prove Proposition~\ref{equality}, for this approach leaves
unaddressed the question of whether inequality~\ref{tapp} alone is restrictive enough to force
rigidity.  When the soul is round, we show that rigidity follows from the inequality alone:
\begin{prop}
Let $f,g:(S^2,\text{round})\ra\R$. Assume that $\int f=0$.  If for all unit tangent vectors $X$,
$(Xf)^2 \leq (f^2 + \frac{2}{3}\text{hess}_g(X,X))$, then (1) equality holds, (2) $f$
is the restriction of a linear function on $\R^3$, and (3) $g$ is uniquely determined
(up to an additive constant) by $f$.
\end{prop}
\begin{proof}
Taking the trace of the inequality gives:
$$|\nabla f|^2\leq 2f^2 + \frac{2}{3}\Delta g.$$
Integrating yields the following:
$$\frac{\int|\nabla f|^2}{\int f^2} \leq 2.$$
The left side of this last inequality is a Rayleigh quotient. Since the lowest non-zero eigenvalue
of $(S^2,\text{round})$ is $2$, and since $f$ is $L^2$-orthogonal to the constant functions, equality must hold, and
$f$ must be an eigenfunction of the laplacian.  That is, $f$ is the restriction to $S^2\subset\R^3$ of a linear function,
so we can write $f(p)=\lb p,Z\rb$ for some vector $Z\in\R^3$.  Notice that $Xf=\lb X,Z\rb$.
If $Z=0$, it is easy to see that $g$ must be constant, so assume $Z\neq 0$.

Let $\gamma:[0,2\pi]\ra S^2$ be a geodesic with $\gamma(0)=\gamma(2\pi)=Z/|Z|$.  Let $g(t)=g(\gamma(t))$ and
$f(t)=f(\gamma(t))$.  Inequality~\ref{tapp} along $\gamma$ becomes:
$$g''(t)\geq f'(t)^2 - f(t)^2.$$
which can be re-written as
$$g''(t)\geq |Z|\cdot(\sin^2(t)-\cos^2(t)) = -|Z|\cos(2t).$$
Since both the left and the right sides of this last inequality integrate to zero, equality must hold.  So,
$$g(t)=\frac{|Z|^2}{4}\cos(2t) + C_1\cdot t + C_2.$$
Since $g$ is periodic, $C_1=0$.  Finally, it is straightforward to see that $C_2$ remains
constant as $\gamma$ varies over all geodesics through $Z/|Z|$.
\end{proof}

Inequality~\ref{tapp} is actually a special case of an inequality which holds in any
dimensions~\cite[Theorem A]{T}.  So it makes sense to ask for metrics of nonnegative curvature
on $S^n\times\R^k$ whether this inequality forces rigidity of the metric at the soul.  The case $n=2,k=3$ is
of particular importance, since it would provide insight into the classification of nonnegatively curved metrics on
$S^2\times S^2$.

\bibliographystyle{amsplain}

\end{document}